\documentclass{article}
\usepackage{amscd}
\usepackage{amssymb}
\newtheorem{lemma}{Lemma}
\newtheorem{theorem}{Theorem}
\newtheorem{corollary}{Corollary}
\newtheorem{definition}{Definition}
\newtheorem{prop}{Proposition}

\newtheorem{example}{Example}

\newcommand{\R}{\mathbb{R}}

\newcommand{\qed}{\nopagebreak\begin{flushright}\framebox[2mm]{}\end{flushright}\vspace{5 mm}}
\title{Mean Curvature Flow, Orbits, Moment Maps}
\author{Tommaso Pacini}
\date{}
\begin{document}
\maketitle
\textbf{Abstract}: Given a compact Riemannian manifold together with a group of isometries, we discuss MCF of the orbits and some applications: eg, finding minimal orbits. We then specialize to Lagrangian orbits in Kaehler manifolds. In particular, in the Kaehler-Einstein case we find a relation between MCF and moment maps which, for example, proves that the minimal Lagrangian orbits are isolated.
\section{Introduction}
If a given submanifold $\Sigma$ in a Riemannian manifold $(M,g)$ is not minimal, 
``mean curvature flow'' (MCF) provides a canonical way to deform it.

Ideally, the flow should exist until either a singularity develops, preventing further 
flow, or the submanifold becomes minimal. In this sense, MCF should be a useful tool 
in the search for minimal submanifolds.

A third, but exceptional, possibility that might occur is exemplified by the 
``translating solitons'' in $\R^n$: submanifolds which, under MCF, simply flow by 
translation and thus never converge to a limiting object.

However, MCF is a difficult topic with many open questions. In particular, there is no general theory which can explain what will happen to all $\Sigma\subseteq (M,g)$ under MCF, or classify which singularities can arise. One is thus forced to study each case individually, or at best to look for classes of submanifolds which, under MCF, behave in the same way. Most often this leads to restrictions on the dimension (eg: curves) or codimension (eg: hypersurfaces) of $\Sigma$, or on properties of the immersion (eg: convex).

In general, the presence of symmetries in a problem reduces the number of variables, 
hopefully making things easier. Regarding MCF, the ``best'' case is when $\Sigma$ is the orbit of a group of isometries of $(M,g)$. It is simple to show that, in this case, all $\Sigma_t$ obtained by MCF are also orbits, and MCF basically reduces to solving an ODE on the (finite-dimensional) space of orbits.

Group actions have been extensively studied. In particular, orbits of a (compact, 
connected) Lie group $G$ acting on a (compact, connected) manifold $M$ can be 
classified into three categories: ``principal'', ``exceptional'' and ``singular''. 
This yields a simple and pretty picture of the geometry of the orbit space $M/G$.

The first goal of this paper is to fit MCF into this framework, analyzing ``what 
happens'' to a principal (or exceptional, or singular) orbit under MCF.

The final picture, presented in theorem \ref{MCF}, constitutes, for several reasons, a 
good ``example zero'' of MCF. It is simple; it generalizes the standard example of the 
``shrinking sphere'' in $\R^n$; it is (co)dimension-independent; and especially, in the 
orbit setting, ``everything we might want to be true for MCF, is true''.

We then restrict our attention to Lagrangian orbits. Using moment maps to 
``get a grasp on them'', we prove that the set $\mathcal{L}(M;G)$ of points
belonging to Lagrangian $G$-orbits constitutes a smooth submanifold in any compact 
Kaehler ambient space. 
When $(M,g)$ is a Kaehler-Einstein (KE) manifold, it was already known that the Lagrangian 
condition is preserved under MCF; we give an independent proof of this in the orbit 
setting, and are then free to apply theorem \ref{MCF} to study how Lagrangian orbits 
evolve under MCF. As a simple corollary, we find that ``backwards MCF'' always leads 
to a minimal Lagrangian orbit.

Futaki proved that compact positive KE manifolds come with a ``canonical'' moment map 
$\mu$. In proposition \ref{prop_formula_H} we show that $\mu$ is intimately related to 
the mean curvature of Lagrangian orbits. As a corollary, we find that minimal Lagrangian orbits (wrt fixed $G$) are isolated.

As already noted, MCF of Lagrangian submanifolds is not a new subject; there is also 
some overlap, in the case of torus actions, with [G]. However, given the number of 
known KE manifolds with large isometry groups, there seems to be no a priori reason 
to limit oneself to tori. Our attempt is to develop a ``complete'' 
picture of the general $G$ case, relying only on the basic tools provided by the 
general theory of $G$-actions, moment maps and transformation groups. In this sense, 
we are not aware of any serious overlap with existent literature.

\textbf{Acknowledgements:} I wish to thank T. Ilmanen for a useful conversation
and P. de Bartolomeis and G. Tian for their long-term support, suggestions and interest. I also gratefully acknowledge the generous support of the University of Pisa and of GNSAGA, and the hospitality of MIT.
\section{Smooth group actions on manifolds} \label{group_actions}
This section is mostly a review of standard facts regarding manifolds with a group action. We refer to [A] for further details.

We adopt the following conventions.
\begin{itemize}
\item $M$ is a compact, connected, smooth manifold. $Diff(M)$ will denote its group of diffeomorphisms.

\item $G$ is a compact, connected, Lie group acting on $M$; ie, we are given a homomorphism $i:G\longrightarrow Diff(M)$. The action of $g\in G$ on $p\in M$ will be denoted $g\cdot p$.
\end{itemize}
The action is ``effective'' if $i$ is injective. Notice that, since $Ker(i)$ is normal in $G$, we may ``reduce'' any $G$-action to an effective $G/Ker(i)$-action.

Whenever a group $H$ is not connected, $H^0$ will denote the connected component containing the identity element.

The action of $G$ on $M$ induces an action of $G$ on $TM$. If $X\in T_pM$, it is defined as follows:
$$g\cdot X:=g_*[p](X)\in T_{gp}M$$
where $g_*$ denotes the differential of the map $g=i(g):M\longrightarrow M$.

Letting $\mathfrak{g}$ denote the Lie algebra of $G$, any $X\in\mathfrak{g}$ induces a ``fundamental vector field'' $\tilde{X}$ on $M$, defined as follows:
$$\tilde{X}(p):=\frac{d}{dt}[exp(tX)\cdot p]_{|t=0}$$
where $exp(tX)$ denotes the 1-dimensional subgroup of $G$ associated to $X$.

For all $p\in M$, we define:

$G\cdot p:=\{g\cdot p:g\in G\}\subseteq M\ \ \ \mbox{``orbit of $p$ (wrt $G$)''}$

$G_p:=\{g\in G:g\cdot p=p\}\leq G\ \ \ \mbox{``stabilizer of $p$ (wrt $G$)''}$

Notice that $G_p$ is a closed subgroup of $G$ and that $G\cdot p\simeq G/G_p$. A generic orbit will often be denoted $\mathcal{O}$. Thanks to the compactness hypotheses, every orbit $\mathcal{O}$ is a smooth embedded submanifold of $M$. 

Notice also that $G_{hp}=h\cdot G_p\cdot h^{-1}$; ie, $G_{hp}$ is conjugate to $G_p$. Thus, to each $p\in M$ we may associate a conjugacy class of subgroups of $G$:
$$\mathcal{O}\mapsto (G_p):=\mbox{ conjugacy class of the stabilizer of any $p\in\mathcal{O}$}$$
This class is called the ``type'' of $\mathcal{O}$.

Let $\mathcal{O}$ be any orbit and $p\in\mathcal{O}$. Let $V:=T_pM/T_p\mathcal{O}$. The action of $G$ on $TM$ restricts to an action of $G_p$ on $T_pM$, and $T_p\mathcal{O}\leq T_pM$ is a $G_p$-invariant subspace. Thus there is a natural action of $G_p$ on $V$.

This induces an action of $G_p$ on $G\times V$, as follows:
$$h\cdot (g,v):=(g\,h^{-1},h\cdot v)$$
Let $G\times_{G_p}V:=(G\times V)/G_p$ denote the quotient space. Then $G\times_{G_p}V$ is a vector bundle (with fiber $V$) over $G/G_p\simeq \mathcal{O}$ and there is an action of $G$ on $G\times_{G_p}V$ as follows:
$$g_1\cdot[g_2,v]:=[g_1\,g_2,v]$$
The following result shows that $G\times_{G_p}V$ contains complete information on the local geometry of the group action near $\mathcal{O}$.
\begin{theorem} \label{theorem_Gnbd} Let $G$, $M$ be as above.  

Then there exist a $G$-invariant neighborhood $U$ of $\mathcal{O}$ in $M$ and a $G$-invariant neighborhood $W$ of the zero section of $G\times_{G_p}V$ such that $U$ is $G$-equivariantly diffeomorphic to $W$.
\end{theorem}
\begin{corollary} \label{cor_orbit_types} Let $M$, $G$ be as above. 
\begin{enumerate} 
\item For each fixed type, the union of all orbits of that type forms a (possibly disconnected) submanifold of $M$.
\item There is only a finite number of orbit types.
\item There is an orbit type $(P)$ whose orbits occupy an open, dense, connected subset of $M$.
\end{enumerate}
\end{corollary}
The types of the $G$-action can be partially ordered by the following relation:
$$\alpha\leq\beta \Leftrightarrow \exists H,K \leq G: \alpha=(H), \beta=(K), H\leq K$$
If a given orbit $\mathcal{O}$ has type $(K)$, any nearby orbit $\mathcal{O}'\subset G\times_KV$ can be written $\mathcal{O}'=G\cdot[1,v]$; it is simple to show that the stabilizer of $[1,v]$ is the stabilizer $K_v$ of $v\in V$ wrt the $K$-action, so it is a subgroup of $K$. In other words, $\mbox{type} (\mathcal{O}')\leq \mbox{type}(\mathcal{O})$.

In particular, the type $(P)$ defined by corollary \ref{cor_orbit_types} must be an absolute minimum:
$$(P)\leq (K), \mbox{   for all types }(K)$$
It is also clear that $\mbox{dim }\mathcal{O}'\geq \mbox{dim }\mathcal{O}$ (the dimension of orbits is a lower-semicontinuous function on $M$) and that orbits of type $(P)$ have maximum dimension among all orbits.

The final picture is thus as follows.

Given $M,G$ as above, there are three categories of orbits:
\begin{enumerate}
\item ``Principal orbits'', corresponding to the minimal type $(P)$.

They occupy an open, dense, connected subset of $M$.
\item ``Exceptional orbits'', corresponding to those types $(K):\ K/P$ is finite. Via the projection $G\times_KV\longrightarrow G/K$, any nearby principal orbit is a finite covering of the exceptional orbit $G/K$.

In particular, exceptional orbits and principal orbits have the same dimension.
\item ``Singular orbits'', corresponding to those types $(K)$: dim $K>$ dim $P$. Their dimension is strictly smaller than that of principal orbits.
\end{enumerate}

\begin{example} \label{example_S2} $S^1$ acts isometrically on $S^2:=\{x^2+y^2+z^2=1\}\subseteq \R^3$ by rotations along the $z$-axis.

The orbits are the sets $S^2\bigcap \pi_c$, where $\pi_c:=\{z\equiv c\}$. The singular orbits, of type ($S^1$), are the poles; all other orbits are principal, of type ($1$).
\end{example}
\begin{example} \label{example_RP2} On $S^2$ there is also an isometric $\mathbb{Z}_2$-action which identifies antipodal points. Since the two actions commute, the $S^1$-action passes to the quotient $\R\mathbb{P}^2\simeq S^2/\mathbb{Z}_2$.

There is one singular orbit, represented by $S^2\bigcap \pi_1$; one exceptional orbit, represented by $S^2\bigcap \pi_0$ (the ``equator''); all other orbits, represented by $S^2\bigcap \pi_c$, $0<c<1$, are principle.
\end{example}
Together, principal and exceptional orbits constitute the set of ``regular orbits''. Any regular orbit $\mathcal{O}=G\cdot q,q\in M^{reg}$ is the image of an immersion 
$$\phi:G/P\hookrightarrow M,\ \ \ [g]\mapsto [g]\cdot q$$
Notice that, if $\mathcal{O}$ is principle, then $\mathcal{O}\simeq G/P$, ie $\phi$ is an embedding. If $\mathcal{O}$ is exceptional, of type $(K)$, then $\mathcal{O}\simeq G/K$ and $\phi$ is a covering map of $G/P$ over $G/K$.

We may set $M^{pr}:=\{p\in M:G\cdot p \mbox{ is a principal orbit} \}\subseteq M$ and, analogously, define $M^{ex},M^{sing}, M^{reg}=M^{pr}\bigcup M^{ex}$.

Each of these subsets, generically denoted $M^*$, is a smooth submanifold inside $M$ and $M^*/G$ also has a smooth structure.

Thus the set $M/G$, which is compact and Hausdorff with respect to the quotient topology, has the structure of a ``stratified smooth manifold'', the smooth strata being the connected components of $M^*/G$. Once again, $M^{pr}/G$ occupies an open, dense, connected subset of $M/G$.
 
An interesting application of all the above is the following, simple, fact.
\begin{corollary} \label{normal_P} Assume $G$ acts on $M$, with principal type $(P)$.
\begin{enumerate}
\item If $P$ is normal, then $P\leq G_p$, $\forall p\in M$.

Thus, if the action is effective, $P=\{1\}$.
\item If $P^0$ is normal, then $P^0\leq G_p$, $\forall p\in M$.

Thus, if the action is effective, $P$ is finite.
\end{enumerate}
In particular, assume a torus $T$ acts effectively on $M$. Then $(P)=\{1\}$.
\end{corollary}
\textbf{Proof}: If $P$ is normal, $P\leq G_p$, $\forall p\in M^{pr}$. Since $M^{pr}$ is dense in $M$, it is easy to prove that $P\leq G_p$, $\forall p\in M$.

The proof of (2) is similar. 
\qed
Our last goal, in this section, is to ``understand'' convergence of orbits.

Assume given a curve of principal orbits $\mathcal{O}_t$ (corresponding to immersions $\phi_t:G/P\hookrightarrow M$) which, in the topology of $M/G$, converges to some limiting orbit $\mathcal{O}$. We must distinguish three cases.

\begin{enumerate}
\item Assume $\mathcal{O}\simeq \phi:G/P\hookrightarrow M$ is principal, ie has minimal type $(P)$. 

Since $M^{pr}$ is open in $M$, each orbit near $\mathcal{O}$ must also have type $(P)$; thus, wrt the local linearization $M=G\times_PV$ based at $\mathcal{O}$, $P$ acts trivially on $V$ and $G\times_PV=G/P\times V$ is the trivial bundle. In particular, this shows that $\mathcal{O}_t\rightarrow \mathcal{O}$ smoothly in $M$, ie $\phi_t\rightarrow \phi$.
\item Assume $\mathcal{O}\simeq\phi$ is exceptional. Then, near $\mathcal{O}$, there are either exceptional or principal orbits and they are coverings of $\mathcal{O}$. It is still true that $\phi_t\rightarrow \phi$ smoothly, but the limit is not injective.
\item Assume $\mathcal{O}$ is singular. Let $K$ be the stabilizer of $p\in\mathcal{O}$.

Locally, $M=G\times_KV,\ \mathcal{O}_t=G\cdot[1,v_t]$ (for some $v_t\in V$) and $K_{v_t}\leq K$ is the stabilizer of $[1,v_t]$. Since $(K_{v_t})\equiv (P)$, all $K_{v_t}$ have constant dimension $q$. The corresponding Lie algebras $\mathfrak{k}_{v_t}$ are thus points in the Grassmannian $Gr(q,\mathfrak{k})$ of $q$-planes in $\mathfrak{k}:=Lie(K)$. By compactness of $Gr(q,\mathfrak{k})$, we may conclude the following: any sequence $\mathcal{O}_n\subseteq\mathcal{O}_t,\  \mathcal{O}_n\rightarrow \mathcal{O}$, contains a subsequence $\mathcal{O}_{n_k}$ such that $\mathfrak{k}_{n_k}\rightarrow \mathfrak{k}_0$, for some $\mathfrak{k}_0\in Gr(q,\mathfrak{k})$.

Let $\{X_1,\dots,X_r\}$ span a complement of $\mathfrak{k}_0$ in $\mathfrak{k}$, and let $\{Y_1,\dots, Y_s\}$ span a complement of $\mathfrak{k}$ in $\mathfrak{g}$. Then $T_p\mathcal{O}$ is generated by the fundamental vector fields $\tilde{Y}_i$ and $T_{[1,v_{n_k}]}\mathcal{O}_{n_k}$ is generated by $\tilde{X}_i,\tilde{Y}_j$. Since $\tilde{X}_i$ are smooth on $M$ and $\tilde{X}_{i|\mathcal{O}}\equiv 0$, we see that $\|\tilde{X}_{|\mathcal{O}_{n_k}}\|\rightarrow 0$ (wrt any invariant metric on $M$).

In other words, convergence to a singular orbit is described, up to subsequences, by the vanishing of certain fundamental vector fields; which fields vanish depends on the particular subsequence.
\end{enumerate}

\section{MCF of orbits}

Let us now fix a compact, connected, Riemannian manifold $(M,g)$ and a compact, connected, Lie group of isometries, $G\leq Isom_g(M)$. $(P)$ will denote the minimal type of the $G$-action, and $\mathfrak{p}$ the corresponding Lie algebra.

Recall that, to any immersion $\phi:\Sigma\hookrightarrow M$, we may associate a volume
$$vol(\phi):=\int_\Sigma vol_{\phi^*g}$$
Since any regular orbit $\mathcal{O}$ corresponds to an immersion $\phi:G/P\hookrightarrow M$, we get a function
$$vol:(M/G)^{reg}\longrightarrow \R, \ \ \ \mathcal{O}\mapsto vol(\phi)$$

The quotient map $\pi: M\longrightarrow M/G$ yields a pull-back map $\pi^*vol:M^{reg}\longrightarrow \R$. We will often write $vol(\mathcal{O})$ instead of $vol(\phi)$ and $vol$ instead of $\pi^*vol$.

\begin{prop} The volume function has the following properties:
\begin{enumerate}
\item $vol: M^{reg} \longrightarrow \R$ is smooth.
\item It has a continuous extension to zero on $M^{sing}$.

This defines a continuous function $vol:M\longrightarrow \R$.
\item The function $vol^2:M\longrightarrow \R$ is smooth.
\end{enumerate}
\end{prop}
\textbf{Proof}: For any regular orbit $\mathcal{O}=\phi:G/P\hookrightarrow M$, $\phi^*g$ defines a $G$-invariant metric on $G/P$. Let $Z_1,\dots,Z_n$ be any basis of $T_{[1]}G/P$, induced by the projection onto $\mathfrak{g}/\mathfrak{p}$ of elements $Z_i\in \mathfrak{g}:Z_i\notin \mathfrak{p}$. Let $\mu:=Z_1^*\wedge\dots\wedge Z_n^*$ be the induced left-invariant volume form on $G/P$. 

Since both volume forms are $G$-invariant, $vol_g=c\cdot \mu$, for some $c=c(\mathcal{O})$. Clearly, $c=\sqrt{det\,g_{ij}}$, where $g_{ij}:=\phi^*g[1](Z_i,Z_j)$. Thus 
$$vol(\mathcal{O})=\int_{G/P}vol_{\phi^*g}=\int_{G/P}\sqrt{det\,g_{ij}}\cdot\mu=\sqrt{det\,g_{ij}}\cdot \int_{G/P}\mu$$
Now let $\mathcal{O}_t$ be a curve of regular orbits, $\mathcal{O}_t=\phi_t:G/P\hookrightarrow M$. Assume that $\mathcal{O}_t\rightarrow \mathcal{O}$. If $\mathcal{O}$ is also regular, $\mathcal{O}=\phi:G/P\hookrightarrow M$, then $p_t:=\phi_t[1]\rightarrow p:=\phi[1]$. Choose $Z_i\in\mathfrak{g}$ such that the induced fundamental vector fields $\tilde{Z}_i$ span $T_p\mathcal{O}$. Then $\tilde{Z}_i$ also span $T_{p_t}\mathcal{O}_t$. Setting $g^t_{ij}:=g[p_t](\tilde{Z}_i,\tilde{Z}_j)$ and $g_{ij}:=g[p](\tilde{Z}_i,\tilde{Z}_j)$, we find 
$$vol(\mathcal{O}_t)=\sqrt{det\,g^t_{ij}}\cdot constant\longrightarrow\sqrt{det\,g_{ij}}\cdot constant\neq 0$$ 
This shows that $vol$ is smooth on $M^{reg}$. If $\mathcal{O}$ is singular, we saw in section \ref{group_actions} that, for any sequence $\mathcal{O}_n\subseteq\mathcal{O}_t:\mathcal{O}_n\rightarrow \mathcal{O}$, we may choose $Z_i$ so that, for some subsequence, certain $\tilde{Z}_i$ vanish. This shows that $\sqrt{det\,g_{ij}^t}\rightarrow 0$, so $vol$ extends continuously to zero on $M^{sing}$.

Since $vol(\mathcal{O}_t)^2=det\,g_{ij}^t\cdot constant$, $vol^2$ is smooth on $M$.
\qed
\begin{corollary} (cfr. [H]) Let $G$ be any compact, connected Lie group acting by isometries on a compact, connected Riemannian manifold $(M,g)$. 

Then there exists a regular minimal orbit of the $G$-action.
\end{corollary}
\textbf{Proof}: Since $M$ is compact, the continuous function $vol:M\longrightarrow \R$ has a maximum, which necessarily corresponds to a minimal (immersed) orbit.
\qed
Example \ref{example_RP2} of section \ref{group_actions} shows that the minimal orbit might be exceptional.

Let us now recall the notion of ``mean curvature flow''.

Fix manifolds $\Sigma$ and $(M,g)$, and an immersion $\phi:\Sigma\hookrightarrow M$.

A smooth 1-parameter family of immersions $\phi_t:\Sigma\hookrightarrow M$ is called a ``solution to the MCF of $(\Sigma,\phi)$'' if it satisfies the following equation

\ 

$\mbox{(MCF)\ \ \ } \left\{ \begin{array}{rcl}
\frac{d}{dt}\phi_t &=& H (\phi_t) \\
\phi_0 &=& \phi
\end{array} \right.$

\ 

where $H(\phi_t)$ denotes the ``mean curvature vector field'' of $\phi_t$, defined as the trace of the second fundamental form of the immersion. It is well-known that $H$ is, up to sign, the ``$L^2$-gradient'' of the volume functional on immersions:
$$\frac{d}{dt}vol(\phi_t)_{|t=0}=-\int_\Sigma(H,\frac{d}{dt}\phi_{t_{|t=0}})$$
Locally, (MCF) can be written as a II-order quasi-linear parabolic system of equations. In particular, solutions always exist for some short time interval $t\in [0,\epsilon)$ and are unique.

We want to focus on solving (MCF) under the assumption that $(\Sigma,\phi)$ is an orbit of a group of isometries.

Consider the map $H:p\mapsto H_p$, which associates to each $p\in M$ the mean curvature $H_p$ of the orbit $G\cdot p$. This defines a vector field on $M$. 

The following lemma examines its continuity/smoothness.

\begin{lemma} Let $H$ be defined as above.
\begin{enumerate} 
\item $H$ is smooth along each submanifold given by orbits of the same type. 

It is also smooth on $M^{reg}$.
\item $H$ is $G$-invariant.
\end{enumerate}
\end{lemma}
\textbf{Proof}: The smoothness of $H$ along orbits of the same type is clear. Smoothness on $M^{reg}$ comes from the convergence properties of regular orbits: basically, $H$ is a local object and does not notice the difference between principal and exceptional orbits. 

(2) is a consequence of the fact that all ingredients in the definition of $H$ are $G$-invariant.
\qed

In particular, $H$ descends to a vector field on $M/G$ and is smooth along each stratum. If $\mathcal{O}=\phi:G/K\hookrightarrow M$ and $p:=\phi([1])$, we can consider the following ODE on $M$:

\ 

$\mbox{(MCF')\ \ \ } \left\{ \begin{array}{rcl}
\dot{p}(t) &=& H[p(t)] \\
p(0) &=& p
\end{array} \right.$

\ 

Notice that, given a solution $p(t)$ of (MCF'), the $G$-equivariant map 
$$\phi:G/K\times [0,\epsilon)\longrightarrow M,\ \ \ \phi([g],t):=g\cdot p(t)$$
solves (MCF) with the initial condition $(\Sigma,\phi)=\mathcal{O}$. By uniqueness of solutions of (MCF), this shows that MCF of an orbit gives a curve of orbits.

In other words, if $(\Sigma,\phi)$ is an orbit, (MCF) is equivalent to the ODE on $M/G$ (or on $M$) determined by integrating $H$.

The reduction of the problem from a PDE to an ODE simplifies things enormously. For example, MCF of orbits has the following properties:
\begin{itemize}
\item There exists a (unique) solution $\mathcal{O}_t$ defined on a maximal time interval $(\alpha,\beta)$: this comes from standard ODE theory.
\item (MCF) may be inverted; ie, $t\mapsto \mathcal{Q}(t):=\mathcal{O}(-t)$ solves the equation for ``backward MCF'':
$$\frac{d}{dt}\mathcal{Q}_t = -H_{\mathcal{Q}_t},\ \ \ \mathcal{Q}_0 = \mathcal{O}$$
This is true for any ODE of the type $\dot{x}=f(x(t))$, but is very atypical for parabolic problems.
\end{itemize}

Another interesting feature of (MCF) on orbits is that it preserves types:
\begin{prop} For each orbit $\mathcal{O}$, $H_\mathcal{O}$ is tangent to the submanifold determined by the type of $\mathcal{O}$.

In particular, if $\mathcal{O}_t$ is the solution of (MCF) with initial condition $\mathcal{O}_0=\mathcal{O}$, then type $(\mathcal{O}_t)\equiv\mbox{type }(\mathcal{O})$.
\end{prop}
\textbf{Proof}: Let $p\in M$ and let $G\cdot p$ have stabilizer $K$. Locally near $G\cdot p$, $M=G\times_KV$, where $V=T_p(G\cdot p)^\perp$ and $K$ acts isometrically on $V$.

This determines a decomposition of $V$ into $K$-irreducible subspaces: $V=\oplus V^i$. Since $H$ is $G$-invariant, it is also $K$-invariant, so $H\in V^0:=\{v\in V:k\cdot v=v,\forall k\in K\}$.

Notice that $G\times_K V^0$ corresponds to the orbits near $G\cdot p$ of type $(K)$. Thus $H_p$ is tangent to the set of such orbits.

Since this is true for each $p\in M$, (MCF) preserves types.
\qed
\begin{corollary} (cfr. [HL]) Let $(M,g)$ be as above. 

If an orbit is isolated wrt all other orbits of the same type, then it is minimal.
\end{corollary}
We now want to show that, on $M^{reg}$, (MCF) is actually a gradient flow; ie, for some $f\in C^\infty(M^{reg})$, $H=\nabla f$.

Let $p_t$ be a curve in $M^{reg}$ and $\mathcal{O}_t:=G\cdot p_t$. We will let $X$ denote both the vector $\frac{d}{dt}p_{t|t=0}$ at $p_0$ and the $G$-invariant vector field $\frac{d}{dt}\mathcal{O}_{t|t=0}$ along $\mathcal{O}$.

Since $H,X$ and the metric on $M$ are $G$-invariant, $(H,X)$ also is. Thus:
$$\frac{d}{dt}vol(p_t)_{|t=0}=\frac{d}{dt}vol(\mathcal{O}_t)_{|t=0}=-\int_{G/P}(H,X)\,vol_\mathcal{O}=-vol(\mathcal{O})\cdot(H,X)_{p_0}$$
This proves that 
$$(H,X)_{p_0}=-\frac{\frac{d}{dt}vol(p_t)_{|t=0}}{vol(p_0)}=-\frac{d}{dt}log\,vol(p_t)_{|t=0}$$
In other words, $H=-\nabla log(vol)$ on $M^{reg}$.

We now have all the information we need to understand how MCF fits into the framework set up in section \ref{group_actions}.

Let $\mathcal{O}=\phi:G/P\hookrightarrow M$ be a fixed principal orbit and let $\mathcal{O}_t=\phi_t:G/P\hookrightarrow M$ be the maximal curve obtained by MCF, with initial condition $\mathcal{O}(0)=\mathcal{O}$ and $t\in(\alpha,\beta)$. Let $p_t=\phi_t([1])$, so that $\frac{d}{dt}p_t=H[p_t]$. 

Since $M$ is compact, there is a sequence $\{p_n\}\subseteq \{p_t\}$ and $\tilde{p}\in M$ such that $p_n\rightarrow \tilde{p}$. Thus $\mathcal{O}_n:=G\cdot p_n\rightarrow \tilde{\mathcal{O}}:=G\cdot \tilde{p}$.

In general, however, different sequences may have different limits, so we cannot hope that $\mathcal{O}_t\rightarrow\tilde{\mathcal{O}}$. The following example of this was suggested to the author by T. Ilmanen.

\begin{example} \label{example_ilmanen}: Consider an embedding $s:\R\hookrightarrow \{(x,y)\in \R^2:x^2+y^2>1\}$ that tends towards $S^1=\{x^2+y^2=1\}\subseteq \R^2$ as $t\rightarrow \infty$, spiralling around it. 

Let $S$ denote its image and $f:S\rightarrow \R^1$ be a positive, decreasing function on $S$ such that $f(t)\searrow c>0$ as $t\rightarrow \infty$.

$S$ may be ``fattened'' by a tubular neighborhood $U$ (of decreasing width, as $t\rightarrow \infty$). At each point $s\in S$, $f$ may be extended, with constant value $f(s)$, in the normal directions. This gives an extension of $f:U\longrightarrow \R$ such that $\nabla f_{|S}$ is tangent to $S$. A partition of unity argument now allows us to extend $f$ to a smooth function $\tilde{f}:\R^2\longrightarrow \R$. Clearly, $\nabla\tilde{f}_{|S}=\nabla f_{|S}$ and $\tilde{f}_{|S^1}\equiv c$.
 
Since we're only interested in what happens near $S^1$, we may perturb $\tilde{f}$ so that it extends to some compact $(\Sigma^2,g)$ containing a neighborhood of $S^1$. We have thus built a smooth gradient vector field on a compact manifold whose flow does not converge to a unique point.

Now let $M:=\Sigma\times S^1$, with the obvious $S^1$-action. For each $p\in \Sigma$, let $\{p\}\times S^1$ have an $S^1$-invariant metric $h_p$ such that $log\,vol(\{p\}\times S^1)=f(p)$.

For each $(p,q)\in M$, let $T_{p,q}M$ have the product metric $g[p]\times h_p[q]$.

Since $H_\mathcal{O}=-\nabla log\,vol(\mathcal{O})=-\nabla f$, MCF of any orbit $\{s\}\times S^1$ with $s=s(t_0)\in S$ yields the curve of orbits $\{s(t)\}\times S^1$, $t\in [t_0,\infty)$. By construction, these orbits have no limit as $t\rightarrow \infty$.
\end{example}
We are thus interested in conditions ensuring the existence of $lim_{t\rightarrow\beta}\mathcal{O}_t$.

The following lemma shows that, in the analytic context, things work nicely.
\begin{lemma} Consider the ODE
$$\dot{p}=-\nabla f,\ \ \ p(0)=p_0$$
with maximal solution $p(t),\ t\in(\alpha,\beta)$.

Assume that, for some subsequence $t_n\rightarrow \beta$, $p(t_n)\rightarrow y$ and that $f$ is analytic in a neighborhood of $y$. Then $p(t)\rightarrow y$.
\end{lemma}
\textbf{Proof}: Let $s(t):=\int_0^t\|\nabla f(p_\tau)\|d\tau$. Since $\|\nabla f(p_\tau)\|>0$, $s$ is a diffeomorphism between $(\alpha,\beta)$ and some $(\alpha',\beta')$. 

Since $f(p_t)$ is monotone, it also gives a diffeomorphism between $(\alpha,\beta)$ and some $(a,b)$, where $b=f(y)$. In particular, $s$ can be written as a function of $f\in(a,b)$; ie, $s=s(f)$.

Since $\frac{df}{dt}=\nabla f\cdot\frac{dp}{dt}=-|\nabla f|^2$, we find $\frac{df}{ds}=-|\nabla f|$ and $\frac{ds}{df}=-\frac{1}{|\nabla f|}$.

It is simple to prove that $p(t)\rightarrow y\Leftrightarrow\beta'<\infty$.

We may assume that $f(y)=0$. When $f$ is analytic near $y$, the ``Lojasiewicz inequality'' asserts that there exists a neighborhood $U$ of $y$ and $\delta>1$, $c>0$ such that, on $U$, $|f|\leq c\,|\nabla f|^\delta$. Thus $\frac{1}{|\nabla f|}\leq c\cdot |f|^{-\frac{1}{\delta}}$, so 
$$\beta'=\int_{f(p_0)}^{f(y)}\dot{s}(f)\,df=\int_b^{f(p_0)}\frac{1}{|\nabla f|}df=\int_0^{f(p_0)}\frac{1}{|\nabla f|}df<\infty$$\qed
If $(M,g)$ is analytic and $vol(\mathcal{O}_t)\geq c>0$, we may apply this lemma to $f:=log\,vol$, proving that $\mathcal{O}^+:=lim_{t\rightarrow\beta}\mathcal{O}(t)$ exists and is regular. We may then apply the following, classical lemma to $N:=M^{reg}$.
\begin{lemma} Let $H$ be a smooth vector field on a manifold $N$. 

Let $p(t):t\in (\alpha,\beta)\longrightarrow N$ be a maximal integral curve and assume that there exists $q\in N$ such that $p(t)\rightarrow q$, as $t\rightarrow \beta$.

Then $\beta=\infty$ and $H(q)=0$.
\end{lemma}
Notice also that $H=-\frac{1}{2}\nabla log(vol^2)=-\frac{1}{2}\frac{\nabla vol^2}{vol^2}$. In the analytic context, this shows that if $\mathcal{O}_t\rightarrow\mathcal{O}$ and $\mathcal{O}$ is singular, then $\|H_{\mathcal{O}_t}\|\rightarrow\infty$ (in the smooth category, following the idea of example \ref{example_ilmanen}, one could build examples for which such a limit does not exist). The following examples show that $\mathcal{O}$ may be minimal or not.
\begin{example}: Let $S^1$ act on $S^2$ as in example \ref{example_S2} of section \ref{group_actions}. Any $\mathcal{O}$ which is not the equator or a pole flows, under MCF, to the closest pole, which is a singular, minimal, orbit. This happens in finite time.
\end{example}
\begin{example}: The above action of $S^1$ on $S^2$ induces an action of $S^1$ on $S^2\times S^2$. Let $\mathcal{O}\simeq S^1\times S^1$ be the product of a ``small'' orbit in $S^2$ (ie: near a pole) and a ``large'' orbit in $S^2$. The flow $\mathcal{O}_t$ becomes singular as soon as the smaller orbit collapses onto the pole, but this limiting curve, $p\times S^1$, is not minimal: its flow exists until the second curve collapses.
\end{example} 
Summarizing, we have proved the following result.
\begin{theorem} \label{MCF} Let $(M,g)$ be a compact, Riemannian manifold and let G be a compact group acting by isometries on $M$. Let $\mathcal{O}$ be a principal orbit. Then:
\begin{enumerate}
\item MCF preserves orbits and types.

Thus there exists a unique, maximal, curve of principal orbits $\mathcal{O}_t$, $t\in (\alpha,\beta)$, solution of MCF with $\mathcal{O}_0=\mathcal{O}$.
\item Assume $(M,g)$ is analytic and $\mathcal{O}^+:=lim_{t\rightarrow \beta} \mathcal{O}_t$ exists. Then:
\begin{itemize} 
\item $\mathcal{O}^+ \mbox{ is a regular orbit } \Leftrightarrow \beta=\infty \Leftrightarrow \|H(t)\|\rightarrow 0$.

In this case, $\mathcal{O}^+$ is minimal.
\item $\mathcal{O}^+ \mbox{ is a singular orbit } \Leftrightarrow \beta<\infty \Leftrightarrow \|H(t)\|\rightarrow \infty$.

In this case, $\mathcal{O}^+$ may be minimal or not.
\end{itemize}
If $vol(\mathcal{O}_t)\geq c>0$, then $\mathcal{O}^+$ always exists and is regular.
\item Assume $(M,g)$ is analytic. Then $\mathcal{O}^-:=lim_{t\rightarrow \alpha} \mathcal{O}_t$ always exists. It is a minimal regular orbit, $\alpha=-\infty$ and $\|H(t)\|\rightarrow 0$. In particular, ``backwards MCF'' always leads to a minimal regular orbit.
\end{enumerate}
\end{theorem}

To get an analogous statement for flows of exceptional or singular orbits, it is sufficient to apply the theorem to the (smooth, compact) manifold $M'$ defined as the closure in $M$ of the set of orbits of the type in question: these orbits will be the principle orbits of the induced $G$-action on $M'$.

\ 

\textbf{Remark}: Using  ``equivariant Morse theory'' applied to the volume function, it would be interesting to study the topology of Riemannian $G$-manifolds in terms of its minimal orbits. In theory, theorem \ref{MCF} would be useful in this.
  
\section{Lagrangian orbits and moment maps} \label{section_moment_maps}

We now want to focus on Lagrangian orbits generated by isometry groups of compact Kaehler manifolds. We start by recalling a few well-known facts concerning transformation groups of Riemannian and Kaehler manifolds. We refer to [K1] for proofs and further details.

\begin{definition} Let $(M,g)$ be a Riemannian manifold. A vector field $X$ on $M$ is an ``infinitesimal isometry'' if $\mathcal{L}_X g\equiv 0$; equivalently, if the local flow generated by $X$ is a curve of isometries.
\end{definition}
$\mathfrak{i}(M)$ will denote the space of all infinitesimal isometries. When $(M,g)$ is complete, $\mathfrak{i}(M)$ is the Lie algebra of $Isom_g(M)$.
\begin{definition} Let $(M,J)$ be a complex manifold. A (real) vector field $X$ on $M$ is an ``infinitesimal automorphism'' if $\mathcal{L}_X J\equiv 0$; equivalently, if the local flow generated by $X$ is a curve of automophisms of $(M,J)$, or if $X-iJX$ is a holomorphic section of $T^{1,0}M$.
\end{definition}
$\mathfrak{h}(M)$ will denote the set of infinitesimal automorphisms. It is closed wrt $J$ and, when $(M,J)$ is compact, it is the complex Lie algebra of the group $Aut_J(M)$ of automorphisms of $(M,J)$.
\begin{theorem} Let $(M,J,g,\omega)$ be a compact Kaehler manifold. Then any infinitesimal isometry is an infinitesimal automorphism, so $Isom_g(M)^0\leq Aut_J(M)^0$.

This, in turn, implies that $Isom_g(M)^0\leq Aut_\omega(M)^0$.
\end{theorem}
The following proposition, although very simple, is the key to understanding Lagrangian orbits.
\begin{prop} \label{prop_lag_orbits} Let $(M^{2n},J,g)$ be a compact Kaehler manifold and let $G\leq Isom_g(M)$ act on $M$ with principal type $(P)$.

Assume there exists a regular Lagrangian $G$-orbit. 

Then $P$ is finite, so dim $G=n$ and $\mathfrak{g}_p=Lie(G_p)=\{0\}$, $\forall p\in M^{reg}$.
\end{prop}
\textbf{Proof}: Assume $\mathcal{O}$ is a Lagrangian orbit. Then $J$ gives an isomorphism $T\mathcal{O}^\perp\simeq T\mathcal{O}\simeq \mathfrak{g}/\mathfrak{p}$ which is equivariant wrt the natural $P$-action. Notice that, for each $p\in P$, this action coincides with the differential of the map
$$p:G/P\longrightarrow G/P,\ \ \ p[g]:=[pg]=[pgp^{-1}]$$
In other words, the action of $p$ on $\mathfrak{g}/\mathfrak{p}$ is the differential of the adjoint action of $p$ on $G/P$. Taken all together, these maps form a group homomorphism $P\longrightarrow GL(\mathfrak{g}/\mathfrak{p})$; the corresponding Lie algebra homomorphism is the map
$$\mathfrak{p}\longrightarrow gl(\mathfrak{g}/\mathfrak{p}),\ \ \ X\mapsto [X,\cdot]$$
If $\mathcal{O}$ is principal, the $P$-action on $T\mathcal{O}^\perp$ is trivial. Thus $P$ acts trivially on $\mathfrak{g}/\mathfrak{p}$ (ie, the action of each $p\in P$ on $\mathfrak{g}/\mathfrak{p}$ is the identity), so the map $\mathfrak{p}\longrightarrow gl(\mathfrak{g}/\mathfrak{p})$ is trivial (ie, the action of each $X\in \mathfrak{p}$ is the zero map), ie $\mathfrak{p}$ is an ideal of $\mathfrak{g}$. This implies that $P^0$ is normal in $G$ and corollary \ref{normal_P} of section \ref{group_actions} proves that $P$ is finite.

Now assume $\mathcal{O}$ is an exceptional Lagrangian orbit of type $(K)$. Locally, $M=G\times_KV$ and $K$ acts as a finite group on $V=\mathfrak{g}/\mathfrak{k}$, so a nbd of $1\in K$ acts trivially on $\mathfrak{g}/\mathfrak{k}$. This shows that $K^0$ acts trivially on $\mathfrak{g}/\mathfrak{k}$.

As above, $K^0$ is normal in $G$. Since $K^0=P^0$, $P^0$ is also normal and we may conclude as above.
\qed
It is now convenient to introduce the concept of Hamiltonian group actions. Again, we refer to [A] for further details.

Let $(M,\omega)$ be a symplectic manifold. Recall that a vector field $X$ on $M$ is ``Hamiltonian'' if $\omega(X,\cdot)$ is an exact 1-form on $M$; ie, $\omega(X,\cdot)=d\,f$, for some $f\in C^\infty(M)$. We say that $f$ is a ``Hamiltonian function'' for $X$.
\begin{definition} The action of $G$ on $M$ is ``Hamiltonian'' if the following conditions are satisfied:
\begin{enumerate}
\item There exists $\mu:M\longrightarrow \mathfrak{g}^*$ such that $<d\,\mu[p](\cdot),X>=\omega[p](\tilde{X},\cdot)$, where $<\cdot,\cdot>$ denotes the natural pairing $\mathfrak{g}^*\times \mathfrak{g}\longrightarrow \R$. Equivalently, $\forall X\in\mathfrak{g},\tilde{X}$ is Hamiltonian (with Hamiltonian function $\mu_X:=p\mapsto <\mu(p),X>$).
\item $\mu$ is $G$-equivariant wrt the $G$-action on $M$ and the coadjoint $G$-action on $\mathfrak{g}^*$. Equivalently, $\mu_{[X,Y]}(p)=\omega[p](\tilde{X},\tilde{Y})$.
\end{enumerate}
We say that $\mu$ is a ``moment map'' for the action.
\end{definition}
\textbf{Remarks}:
\begin{enumerate}
\item Assume $(M,J,g,\omega)$ is a compact Kaehler manifold and that, for some $G\leq Isom_g(M)^0$, condition (1) above is satisfied. Then, $\forall X\in \mathfrak{g}$, 
$$d\,\mu_X=\omega(\tilde{X},\cdot)=g(J\tilde{X},\cdot)$$
This shows that $\nabla \mu_X=J\tilde{X}$, so $\nabla \mu_X$ is an infinitesimal automorphism of $(M,J)$.
\item By definition, the differential $d\,\mu[p]:T_pM\longrightarrow \mathfrak{g}^*$ is the dual of the map
$$d\,\mu[p]^*:\mathfrak{g}\longrightarrow (T_pM)^*,\ \ \ X\longrightarrow d\,\mu_X[p]$$
Thus $Im\,d\,\mu[p]=(Ker\ d\,\mu[p]^*)^\#=(\mathfrak{g}_p)^\#$, where $\mathfrak{g}_p=Lie(G_p)$.

In particular, $d\,\mu[p]$ is surjective iff $G_p$ is discrete.
\end{enumerate}
\begin{definition} $\Sigma\subseteq (M,\omega)$ is ``isotropic'' if $\omega_{|\Sigma}\equiv 0$; if dim $\Sigma=n$ and dim $M=2n$, then isotropic submanifolds are called ``Lagrangian''.
\end{definition}
We are mainly interested in moment maps for the following reason.
\begin{lemma} \label{isotropic_orbits} Assume the action of $G$ on $(M,\omega)$ is Hamiltonian, with moment map $\mu$. Let $p\in M$. Then the following conditions are equivalent:
\begin{enumerate}
\item $\mu$ is constant on the orbit $\mathcal{O}=G\cdot p$.
\item $\mathcal{O}$ is isotropic. 
\item $\mu(p)\in[\mathfrak{g},\mathfrak{g}]^\#$.
\end{enumerate}
\end{lemma}
We now have all the elements necessary to prove the following
\begin{corollary} Let $(M,J,g)$ be a compact Kaehler manifold. Assume $G\leq Isom_g(M)$ acts in a Hamiltonian fashion. Then the set
$$\mathcal{L}(M;G):=\{p\in M^{reg}:G\cdot p\mbox{ is a Lagrangian orbit}\}$$
either is empty or is a smooth submanifold of $M^{reg}$, of dimension $2n-\mbox{dim }[\mathfrak{g},\mathfrak{g}]$.
\end{corollary}
\textbf{Proof}: If there exists a regular Lagrangian orbit, then, by proposition \ref{prop_lag_orbits}, $P$ is finite and dim $G$=n.  Thus every regular isotropic orbit has dimension $n$ and is Lagrangian. Lemma \ref{isotropic_orbits} now shows that, if we let $\mu_{reg}$ denote the restriction of $\mu$ to $M^{reg}$, then $\mathcal{L}(M;G)=\mu_{reg}^{-1}([\mathfrak{g},\mathfrak{g}]^\#)$. Since $\mathfrak{g}_p=0$, $\mu_{reg}$ is a submersion so $\mathcal{L}(M,G)$ is smooth, of dimension $n+\mbox{dim }[\mathfrak{g},\mathfrak{g}]^\#$.
\qed
\begin{example} Assume $G\leq Isom_g(M)$ is semisimple. Then the $G$-action on $M$ is Hamiltonian (cfr. [A]) and Lagrangian orbits are isolated.
\end{example}
\begin{example} \label{ex_torus} Assume that a torus $T^n\leq Isom_g(M)$ acts effectively on $M$ and that $H^1(M;\R)=0$. Then the action is Hamiltonian (cfr. [A]), $P={1}$, $[\mathfrak{g},\mathfrak{g}]=0$ and every regular orbit is Lagrangian. In other words, $\mathcal{L}(M,G)=M^{reg}$.

In particular, there exists a minimal, Lagrangian orbit (cfr. also [G]).

An example of this is provided by the standard $T^n$-action on $\mathbb{P}^n$.
\end{example}

\section{MCF of Lagrangian orbits in KE manifolds} 

In this section, we will assume that $(M,J,g,\omega)$ is a KE manifold.

Since we are interested in group actions, we must recall (cfr. [K1]) some basic facts concerning their transformation groups.

\begin{theorem} \label{KE+} Let $M$ be a compact KE manifold such that $Ric=c\cdot g$, $c>0$.

For any (real) vector field $X$, let $Z:=X-iJX$ and $\zeta:=g(Z,\cdot)$. Then
\begin{enumerate}
\item $\mathfrak{i}(M)$ is totally real in $\mathfrak{h}(M)$: ie, if $X\in\mathfrak{i}(M)$, then $JX\notin \mathfrak{i}(M)$.
\item $X\in \mathfrak{h}(M)\Leftrightarrow \zeta=\overline{\partial}f:f\in C^\infty(M;\mathbb{C}), \Delta f=2cf$.

In particular, $\int_M f=0$, so such an $f$ is unique.
\item $X\in\mathfrak{i}(M)\Leftrightarrow Re(f)=0$.

If we set $E_{2c}:=\{f\in C^\infty(M;\mathbb{R}):\Delta f=2cf\}$, there is an isomorphism:
$$E_{2c}\simeq\mathfrak{i}(M),\ \ \ f\mapsto i\overline{\partial}f=\zeta$$
\item $\mathfrak{h}(M)=\mathfrak{i}(M)\oplus J\,\mathfrak{i}(M)$.
\end{enumerate}
\end{theorem}
It is possible (cfr. [K2]) to prove that positive compact KE manifolds are simply connected. This implies that every fundamental vector field induced by $G\leq Isom_g(M)$ is Hamiltonian. From our point of view, however, much more is true:
\begin{prop} (cfr. [F]) \label{fukaya} Let $M$ be a compact positive KE manifold and $G\leq Isom_g(M)$. Then the action of $G$ on $M$ is Hamiltonian.

Recall the correspondence and the notation from theorem \ref{KE+} above:
$$X\in\mathfrak{i}(M)\leftrightarrow f:f\in C^\infty(M;\R), \Delta f=2cf, \zeta=i\overline{\partial}f$$
Then $\mu_X:=-\frac{1}{2}f$ defines a moment map for $M$, $G$.
\end{prop}
Moment maps are usually not unique: if $\mu$ is a moment map, $\mu+c$ also is, for any $c\in[\mathfrak{g},\mathfrak{g}]^\#\leq \mathfrak{g}^*$. The proposition above suggests the following
\begin{definition} The moment map defined in proposition \ref{fukaya} above will be called the ``canonical moment map'' of the $G$-action.
\end{definition}
Recall, however, that moment maps are uniquely defined on $[\mathfrak{g},\mathfrak{g}]$ because $\mu_{[X,Y]}=\omega(\tilde{X},\tilde{Y})$. Recall also lemma \ref{isotropic_orbits}. Proposition \ref{fukaya} thus leads to the following result:
\begin{corollary} Let $M$ be a compact KE manifold such that $Ric=c\cdot g, c>0$. 
\begin{enumerate}
\item $\forall X,Y\in \mathfrak{i}(M)$, $\omega(\tilde{X},\tilde{Y})\in E_{2c}$.
\item $\forall f\in E_{2c}$, $f$ restricted to $\mathcal{L}(M;G)$ is $G$-invariant.
\end{enumerate}
\end{corollary}
Putting everything together and using the fact that KE metrics are analytic, we can now prove the following result:
\begin{theorem} \label{Lag_MCF} Let $M$ be a compact positive KE manifold and let $G\leq Isom_g(M)$. 

Assume $\mathcal{L}(M;G)$ is not empty. Then $H$ is tangent to $\mathcal{L}(M;G)$, so MCF preserves the Lagrangian condition and may be studied as in theorem \ref{MCF}. 

Furthermore, there exists a minimal Lagrangian orbit.
\end{theorem}
\textbf{Proof}: Recall, for any Lagrangian submanifold $\Sigma$ immersed in Kaehler
$M$, the isomorphism
$$(T\Sigma)^\perp\simeq \Lambda^1(\Sigma),\ \ \ V\simeq\nu:=\omega(V,\cdot)_{|\Sigma}$$
It is well-known (cfr. [TY]) that if $\sigma_H\in \Lambda^1(\Sigma)$
denotes the 1-form corresponding to the mean curvature vector field $H$ under 
this isomorphism, then $d\,\sigma_H=\rho_{|\Sigma}$, where $\rho(X,Y):=Ric(JX,Y)$
is the ``Ricci 2-form" of $M$. 

When $M$ is KE, $\rho=c\cdot\omega$, so this shows that $\sigma_H$ is closed.

Now let $p\in\mathcal{L}(M;G)$. Then 
$T_p\mathcal{L}=\{X\in T_pM:d\,\mu[p](X)\in[\mathfrak{g},\mathfrak{g}]^\#\}$,
so we need to prove that $d\,\mu[p](H)\in[\mathfrak{g},\mathfrak{g}]^\#$.

Since $H$ is $G$-invariant, $\sigma_H$ also is; ie, $\sigma_H\in\mathfrak{g}^*$.
Notice that $d\,\mu[p](H)=\omega[p](\cdot,H)=-\sigma_H[p]$. 

Recall that, for any 1-form $\alpha\in\Lambda^1(\Sigma)$,
$$d\,\alpha(X,Y)=X\alpha(Y)-Y\alpha(X)-\alpha[X,Y]$$
Thus $0=d\,\sigma_H(X,Y)=-\sigma_H[X,Y],\ \forall X,Y\in\mathfrak{g}$,
as desired. 

(1) is now obvious. The properties of $vol$ show that there is a Lagrangian orbit $\mathcal{O}$ of maximum volume (which is minimal in $\mathcal{L}(M;G)$). Let $\mathcal{O}_t$ be a curve in $\mathcal{L}(M;G)$ such that $\mathcal{O}_0=\mathcal{O}$ and $\frac{d}{dt}\mathcal{O}_t=H$. Then $0=\frac{d}{dt}vol(\mathcal{O}_t)_{|t=0}=-\int_M(H_\mathcal{O},H_\mathcal{O})$, so $H_\mathcal{O}\equiv 0$.
\qed
\textbf{Remark}: When $M$ is compact Kaehler Ricci-flat, one can show that $Isom_g(M)^0$ is a torus, so example \ref{ex_torus} show that the analogous statement is trivially true. 

When $M$ is compact negative KE, $Isom_g(M)^0=\{Id\}$, so these manifolds are not interesting from our point of view. Cfr. [K1] for details.

\ 

Our final goal is to explore the relationship between MCF and the canonical moment map.
\begin{prop} \label{prop_formula_H} Let $M^{2n}$ be a compact KE manifold such that $Ric=c\cdot g,c>0$. Given $G\leq Isom_g(M)$, assume $\mathcal{L}(M;G)$ is not empty. Let $\mu:M\longrightarrow \mathfrak{g}^*$ denote the canonical moment map. Then, on $\mathcal{L}(M;G)$, 
\begin{enumerate}
\item $\forall X\in\mathfrak{g},\ \ \ H_\mathcal{O}\cdot\nabla\mu_X=c\,\mu_X$.
\item $\forall p\in \mathcal{L}(M;G)$, the natural ($G$-invariant) metric on $G\cdot p\subseteq M$ defines metrics on $\mathfrak{g}$ and $\mathfrak{g}^*$. For the induced norm (which depends on $p$),
$$d\|\mu\|^2[p](H)=2c\|\mu(p)\|^2$$
\end{enumerate}
\end{prop}
\textbf{Proof}: Let $\mathcal{O}\simeq G\cdot p/G_p$ denote any regular Lagrangian orbit. Let $e_1,\dots,e_n\in\mathfrak{g}\simeq T_p\mathcal{O}$ be a orthonormal basis wrt the induced metric. To simplify the notation, we will denote the corresponding fundamental vector fields also by $e_i$. Then
$$(H_\mathcal{O},\nabla\mu_X)=(\nabla^\perp_{e_j}e_j,\nabla\mu_X)=(\nabla_{e_j}e_j,\nabla\mu_X)=-(e_j,\nabla_{e_j}\nabla\mu_X)$$
In section \ref{section_moment_maps}, we saw that $\nabla\mu_X$ is an infinitesimal automorphism of $M$. Thus
$$\nabla_{Je_j}\nabla\mu_X=\nabla_{\nabla\mu_X}Je_j+[Je_j,\nabla\mu_X]=J(\nabla_{\nabla\mu_X}e_j+[e_j,\nabla\mu_X])=J(\nabla_{e_j}\nabla\mu_X)$$
The definition of the canonical moment map now shows that
\begin{eqnarray*}
2(H_\mathcal{O},\nabla\mu_X)&=&-(e_j,\nabla_{e_j}\nabla\mu_X)-(Je_j,J(\nabla_{e_j}\nabla\mu_X))\\
&=&-(e_j,\nabla_{e_j}\nabla\mu_X)-(Je_j,\nabla_{Je_j}\nabla\mu_X)\\
&=&-div_M(\nabla\mu_X)=\ \Delta_M\mu_X=\ 2c\,\mu_X
\end{eqnarray*}
This proves (1). Applying (1) to $X=e_i$, multiplying by $2\mu_{e_i}$ and summing wrt $i$ shows that
$$H_\mathcal{O}\cdot\nabla\|\mu\|^2=2c\|\mu\|^2$$
which is (2).\qed
We can now prove
\begin{theorem} \label{main}
Let $M^{2n}$ be a compact KE manifold such that $Ric=c\cdot g,c>0$. For $G\leq Isom_g(M)$, let $\mu$ denote the canonical moment map and let $E_{2c}(G):=\{f\in E_{2c}:f=\mu_X,\mbox{ for some }X\in\mathfrak{g}\}$.

Assume that regular orbits have dimension $n$. Then a Lagrangian orbit $\mathcal{O}$ is minimal iff $\mu(\mathcal{O})=0$. In particular, minimal Lagrangian orbits are isolated. Furthermore, the following are equivalent:
\begin{enumerate}
\item There exists a Lagrangian orbit.
\item There exists a minimal Lagrangian orbit.
\item $0\in \mu(M)$.
\item The set $\{p\in M: f(p)=0, \ \forall f\in E_{2c}(G)\}$ is not empty.
\end{enumerate}
\end{theorem}
\textbf{Proof}: By hypothesis, an orbit is regular iff it is $n$-dimensional. In particular, every Lagrangian orbit $\mathcal{O}$ is regular. We may thus restrict our attention to $M^{reg}$. 

If $\mathcal{O}$ is minimal, proposition \ref{prop_formula_H} shows that $\mu(\mathcal{O})=0$. Viceversa, assume that $\mu(\mathcal{O})=0$. Let $\mathcal{O}_t$ be obtained by MCF applied to $\mathcal{O}$. Then proposition \ref{prop_formula_H} shows that $f(t):=\|\mu\|^2(\mathcal{O}_t)$ satisfies
$$\frac{d}{dt}f(t)=2cf,\ \ \ f(0)=0$$
This implies that $f(t)\equiv 0$, so $\mathcal{O}(t)\subseteq \mu^{-1}(0)$. However, $\mu$ is a submersion, so $\mu^{-1}(0)$ is smooth of dimension $n$ and, since $P$ is finite, the elements of $\mu^{-1}(0)/G$ are isolated. Thus $\mathcal{O}(t)\equiv\mathcal{O}$, ie $\mathcal{O}$ is minimal.

Together with theorem \ref{Lag_MCF}, this shows that (1), (2) and (3) are equivalent. The equivalence of (3), (4) comes directly from the definition of $\mu$.
\qed
\textbf{Remark}: In the toric case, one can show that $\mu^{-1}(0)$ is connected, so theorem \ref{main} implies that the minimal Lagrangian orbit is unique. This result was obtained also in [G], by lifting the $T^n$-action from $M$ to its canonical bundle $K_M$ and studying the induced geometry.

\ 

\ 

\mbox{\Large\textbf{Bibliography}}
\begin{description}
\item[\mbox{[A]}] Audin, M., The topology of torus actions on symplectic manifolds, Birkh\"{a}user, 1991
\item[\mbox{[F]}] Futaki, A., The Ricci curvature of symplectic quotients of Fano manifolds, Tohoku Math. J., 39 (1987), 329-339
\item[\mbox{[G]}] Goldstein, E., Minimal Lagrangian tori in Kaehler Einstein manifolds, math.DG/0007135 (pre-print)
\item[\mbox{[H]}] Hsiang, W., On the compact homogeneous minimal submanifolds, Proc. Nat. Acad. Sci., 56 (1966), pp. 5-6
\item[\mbox{[HL]}] Hsiang, W. and Lawson, H.B. Jr., Minimal submanifolds of low cohomogeneity, J. Diff. Geom., 5 (1971), pp. 1-38 
\item[\mbox{[K1]}] Kobayashi, S., Transformation groups in differential geometry, Springer, 1972
\item[\mbox{[K2]}] Kobayashi, S., On compact Kaehler manifolds with positive definite Ricci tensor, Ann. of Math. (2), 74 (1961), pp. 570-574
\item[\mbox{[TY]}] Thomas, R. and Yau, S.-T., Special Lagrangians, stable bundles and mean curvature flow, math.DG/0104197 (preprint)
\end{description} 

\ 

\ 

Tommaso Pacini (Imperial College/University of Pisa)

\ 

Email: pacini@paley.dm.unipi.it

\ 

Subject Class: 53C44 (primary), 53D20 (secondary)
\end{document}